\documentclass[11pt,a4paper]{amsart}

\usepackage{graphicx}

\usepackage{amsmath}%
\usepackage{amssymb,indentfirst,calc,mathrsfs,euscript,bbm}%

\usepackage{setspace}%
\usepackage[reals]{layout}%
\usepackage{xr}%
\usepackage{amscd}%
\usepackage{comment}
\usepackage{epstopdf}
\usepackage{tikz}

\def\square{\hfill\hbox{\vrule height .9ex width .8ex depth -.1ex}}

\def\nit{\hbox{\it I\hskip -2pt  N}}

\newtheorem{lem}{Lemma}[section]

\newtheorem{theo}{Theorem}[section]

\begin{document}

\title[AS]{ Asymptotic properties of optimal trajectories in dynamic programming}
\author{
Sylvain Sorin, Xavier Venel, Guillaume Vigeral }
\address{
Equipe Combinatoire et Optimisation, CNRS FRE 3232,  Facult\'e de Math\'ematiques, UPMC-Paris 6,
  175 Rue du Chevaleret,
75013 Paris, France \newline
GREMAQ
Université de Toulouse 1
Manufacture des Tabacs, Aile J.J. Laffont
21 allée de Brienne
31000 Toulouse, France \newline
INRIA Saclay - Ile-de-France  and  CMAP, Ecole Polytechnique, route de Saclay, 91128 Palaiseau cedex, France
 } \email{sorin@math.jussieu.fr, xavier.venel@sip.univ-tlse1.fr,
guillaumevigeral@gmail.com}
\date{ October 2009}

\bibliographystyle{apalike}

\begin{abstract}

We show in a dynamic programming framework  that uniform convergence of the finite horizon values implies that  asymptotically the average accumulated payoff is constant on optimal trajectories. We analyze and discuss several possible extensions to two-person games.
\end{abstract}

\maketitle

\section{Presentation}
Consider a dynamic programming problem  as described in Lehrer and Sorin  \cite{LS}. Given a set of states $S$,  a correspondence $\Phi$ from $S$ to itself with non empty values   and a payoff function $f$ from $S $ to $[0,1]$,  a feasible play  at $s\in S$ is a sequence $\{ s_m\}$ of states with $s_1= s$ and $s_{m+1} \in \Phi (s_m) $. It induces  a sequence of payoffs $\{ f_m = f(s_m)\}, m = 1, ..., n , ...$.
Recall that starting from a standard problem with random transitions and/or signals on the state, this presentation amounts to work on the set of probabilities on  $S$ and to consider  expected payoffs.\\

Let $v_n(s)$ (resp. $v_{\lambda}(s)$)  be the value of the $n$ stage program  $G_n(s)$  (resp. $\lambda$ discounted program $G_{\lambda}(s)$) starting from state $s$. The  {\bf asymptotic approach}
deals with asymptotic properties of the  values  $v_n$ and $v_{\lambda}$  as $n$ goes to $\infty$ or $\lambda$ goes to  0.\\

The {\bf uniform approach} focuses on properties of the strategies that hold uniformly in long horizons.
$v_{\infty}$ is the uniform value if  for  each $\varepsilon >0$ there exists $N$  such that  for each $s\in S$:\\
1) there is a feasible play $\{ s_m\}$ at $s$ with
$$
\frac{1}{n} \sum_{m=1}^n f(s_m) \geq v_{\infty}(s) - \varepsilon,  \qquad \forall n \geq N
$$
2) for any feasible play  $\{ s'_m\}$ at $s$ and any  $n \geq N$
$$
\frac{1}{n} \sum_{m=1}^n f(s'_m) \leq v_{\infty}(s) + \varepsilon.
$$

Obviously the second  approach is more powerful than the second  (existence of a uniform value implies existence of an asymptotic value : the limit of $v_n$ exists) but it is  also more demanding: there are problems  without uniform value where the asymptotic value exists (see Section 2). \\
Note that the condition for the existence of a uniform value implies that the average accumulated payoff  on optimal trajectories remains close to the value.

We will prove that a similar phenomenon holds true under conditions that are stronger than
 the existence of an asymptotic value but weaker than the existence of a uniform value.

Say that the dynamic programming problem is {\bf regular} if :\\
i) $\lim v_n (s) = v(s)$ exists for each $s\in S$.\\
ii) the convergence is uniform.\\
This condition was already introduced and studied in Lehrer and Sorin \cite{LS} (see Section 2).

We  consider  the following property  {\bf P}:\\

For any $ \varepsilon > 0$, there exists $n_0$,  such that for all $n\geq n_0$,  for any state $s$ and any feasible play $\{ s_m\}$ $\varepsilon$-optimal for $G_n (s)$  and  for any $t\in [0,1]$:
\begin{equation}\label{1}
3\varepsilon\geq  \frac {1}{n} (  \sum_{m=1}^{[tn]} f_m ) - tv(s) \geq -3 \varepsilon.
\end{equation}
where $[tn]$ stands for the integer part of $tn$.\\

This condition says that the average payoff remains close to the value on every almost-optimal trajectory with long duration (but the trajectory may depend on this duration).\\
It also implies a similar property on every time interval.\\

\section{Examples and comments}
1) The existence of the asymptotic value  is not enough to control the payoff as required in property  {\bf P}.\\
 An example  is given in Lehrer and  Sorin \cite{LS} (Section 2),  where both $\lim v_n$ and $ \lim v_{\lambda}$ exist on $S$ but where  the asymptotic average payoff is not constant on the unique optimal trajectory, nor on $\varepsilon$-optimal  trajectories: in  $G_{2n}$,  an optimal play will induce $n$  times 0 then $n$ times 1 while $v = 1/2$.\\
 Note that this example is not regular: the convergence of $v_n$ to $v$  is not uniform.

2) Recall  that in the framework of dynamic programming, regularity  is also equivalent to  uniform convergence of $v_{\lambda}$ (and with the same limit), see Lehrer and  Sorin \cite{LS} (Section 3).\\
Note also that this regularity condition is not sufficient to obtain the existence of a uniform value, see Monderer and  Sorin \cite{MS} (Section 2).

3) General conditions for regularity can be found in Renault \cite{JR}.

\section{Main result}

\begin{theo}\label{a}
Assume that the program is regular, then {\bf P} holds.
\end{theo}

{\bf Proof}

Let us start with the upper bound inequality  in  (\ref{1}). \\
The result is clear for $t \leq \varepsilon$ (recall that that the payoff is in $[0,1]$). Otherwise let $n_1$ large enough so that $n\geq n_1$ implies  $||v_{n}- v || \leq  \varepsilon$ by uniform convergence. Then the required inequality holds  for $ n \geq n_2$ with $[\varepsilon n_2] \geq n_1$.

Consider now the lower bound inequality in  (\ref{1}). \\
The result holds for $t \geq 1-\varepsilon$ by the $\varepsilon$-optimal property of  the play, for $n\geq n_1$.
Otherwise  we use the following lemma from Lehrer and Sorin \cite{LS} (Proposition 1).
\begin{lem}
Both $\limsup v_n$ and $\limsup v_{\lambda}$ decresase on feasible histories.
\end{lem}
In particular, starting from $s_{[tn]}$ the value of the program for the last $n - [t n]$ stages is at most $v (s_{[tn]}) + \varepsilon$  for $n\geq n_2$, by uniform convergence,  hence less than the initial $v(s) + \varepsilon$, using the previous Lemma.
Since the play is $\varepsilon$-optimal in $G_n(s)$, this implies that
\begin{equation}\label{2}
\sum _{m= 1}^{[tn]} f_m + (n - [tn]) (v(s) + \varepsilon )  \geq n(v_n(s) - \varepsilon) \geq   n(  v(s) - 2 \varepsilon)
\end{equation}
hence the required inequality.
\square

\section{Extensions}

\subsection{Discounted case} {} \\
A similar result  holds for the program $G_{\lambda}$ corresponding to the evaluation $\sum_{m=1}^{\infty} \lambda (1- \lambda)^{m-1}f_m$. Explicitly,  one introduces the property ${\bf P'}$:\\
For any $ \varepsilon > 0$, there exists $\lambda_0$,  such that for all $\lambda \leq
\lambda_0$,  for any state $s$ and any feasible play $\{ s_m\}$ $\varepsilon$-optimal for $G_{\lambda} (s)$  and  for any $t\in [0,1]$:
\begin{equation}\label{3}
%{\bf E}_{\sigma_n, \tau_n} ^{\omega}
3\varepsilon\geq    \sum_{m=1}^{n(t;\lambda)}\lambda (1- \lambda)^{m-1} f_m ) - tv(s) \geq -3 \varepsilon.
\end{equation}
where $n(t;\lambda) = \inf \{p \in  \nit; \sum_{m=1} ^p\lambda (1- \lambda)^{m-1} \geq t \}$. Stage $n(t;\lambda)$ corresponds to the fraction $t$ of the total duration of the program.
\begin{theo}
Assume that the program is regular, then {\bf P'} holds.
\end{theo}

{\bf Proof}

The proof follows the same lines than the proof of Theorem \ref{a}.\\
Recall that by regularity both $v_n $ and $v_{\lambda}$ converge uniformly to $v$.
Moreover the discounted sums   $ (1-\lambda)^{-N}\sum _{m=1} ^N \lambda (1- \lambda)^{m-1} f_m$ belong to the convex hull of the averages $\frac{1}{n} \sum _{m=1} ^n f_m; 1\leq n\leq N$.
The counterpart of equation (\ref{2}) is now
\begin{equation}\label{4}
 \sum_{m=1}^{n(t;\lambda)}\lambda (1- \lambda)^{m-1} f_m  + (1-t) (v(s) + \varepsilon )  \geq (v_{\lambda}(s) - \varepsilon) \geq     v(s) - 2 \varepsilon
\end{equation}
\square

\subsection{Continuous time}
{ }  \\
Similar results holds in the following set-up:
$v_T(x)$ is the value of the control problem $\Gamma_T$ with control set $U$  where the state variable in $X$   is governed by a differential equation (or more generally a differential inclusion)
$$
\dot x_t = f(x_t, u_t)
$$
starting from $x$ at time 0. The real  payoff function is $g(x,u)$ and the evaluation is given by:
$$
\frac{1}{T}
\int_0^T g(x_t, u_t) dt.
$$
Regularity in this framework amounts to uniform convergence  (on $X$) of  $V_T$ to some $V$. (Sufficient conditions for regularity can be found in Quincampoix and Renault \cite{QR}).\
The corresponding property is now {\bf P"}:\\
For any $ \varepsilon > 0$, there exists $T_0$,  such that for all $T\geq T_0$,  for any state $x$ and any feasible trajectory $\varepsilon$-optimal for $\Gamma_ T  (x)$  and  for any $\theta \in [0,1]$:
\begin{equation}\label{cc}
3\varepsilon\geq  \frac{1}{T}
\int_0^{\theta T} g(x_t, u_t) dt - \theta V(x) \geq -3 \varepsilon.
\end{equation}

\begin{theo}
Assume that the optimal control problem  is regular, then {\bf P"} holds.
\end{theo}
{\bf Proof}

Follow exactly the same lines than the proof of Theorem (\ref{2}).
\square
\\

Finally the same tools can be used for an evaluation of the form $\lambda \int_0^{+\infty} e^{- \lambda t} g(x_t, u_t) dt$.

\section{Two-player zero-sum games}
In trying to extend this result to a two-person zero-sum framework, several problems occurs.

\subsection{Optimal strategies on both sides}  { } \\
First it is necessary,  to obtain good properties on the trajectory, to ask for optimality on both sides.\\
For example in the Big Match with no signals,
$$
\begin {tabular}{cccc}
& \multicolumn{1}{c}{$$} &\multicolumn{1}{c}{$\alpha$} &\multicolumn{1}{c}{$\beta$} \\
\cline{3-4}
& \multicolumn{1}{c}{$a$} &\multicolumn{1}{|c}{$1^*$} &
\multicolumn{1}{|c|}{$0^*$} \\
\cline{3-4}
& \multicolumn{1}{c}{b} &\multicolumn{1}{|c}{$0$} &
\multicolumn{1}{|c|}{$ 1$}\\
\cline{3-4}
\end{tabular}
$$
where a $*$ denotes an absorbing payoff,
the  optimal strategy of player 1 in the ``asymptotic game" on $[0,1]$ is to play ``$a$ before time $t$" with probability   $t$, see Sorin \cite{S} Section 5.3.2. Obviously, if there is no restrictions on player 2's moves the average payoff will not be constant. However, the optimal strategy  of player 2 is ``always $(1/2,1/2)$"  hence time independent  on $[0,1]$. It thus  induces  a constant payoff and it is easy to see that the property is robust to small perturbations in the evaluation of the payoff.

\subsection{Player 1 controls the transition.} { } \\
Consider a repeated game with finite characteristics (states, moves, signals, ...) and use the recursive formula corresponding to the canonical representation with entrance laws being consistent probabilities on the universal belief space, see Mertens, Sorin and Zamir \cite{MSZ}, Chapters III.1, IV.3. This representation preserves the values but in the auxiliary game, if player 1 controls the transition an optimal strategy of player 2 is to play a stage by stage best reply. Hence the model reduces to the dynamic programming framework and the results of the previous sections apply.\\
A simple example corresponds to a game with incomplete information on one side where asymptotically an optimal strategy of the uniform player 1 is a splitting at time 0, while player 2 can obain $u(p_t)$ at time $t$ where $u$ is the value of the non-revealing game and $p_t$ the martingale of posteriors at time $t$, see Sorin \cite{S}, 3.7.2.

\subsection{Example.} { } \\
Back to the general framework of two person zero-sum repeated games, the following example shows that in addition one has to strengthen the conditions on the pair of $\varepsilon$-optimal strategies. We exhibit a game having a uniform value $v$ but  for some state $s$ with $v(s)=0$ one can construct, for each $n$,  optimal strategies in $\Gamma_n(s)$ inducing essentially a  constant payoff 1 during the first half of the game.\\

Starting from the initial state $s$,  the tree representing the game $ \Gamma$ has countably many subgames $\tilde\Gamma_{2n}$, the
transition being controlled by player 1 (with payoff 0). In
$\tilde\Gamma_{2n}$ there are at most $n$ stages before reaching an
absorbing state.  At each of these stages of the form $(2n,m), m = 1, ...n,$ the players plays
a ``jointly controlled" process  leading  either to a payoff 1 and
the next stage $(2n,m+1)$  (if they agree) or an absorbing payoff
$x_{2n,m}$ with  $(m-1) + (2n -(m-1))x_{2n,m} = 0$, otherwise.  Hence  every feasible path of length $2n$ in $\tilde
\Gamma_{2n}$ gives a total payoff 0. Obviously the uniform value
exists since each player can stop the game at each node, inducing
the same absorbing payoff. The representation is as follows:

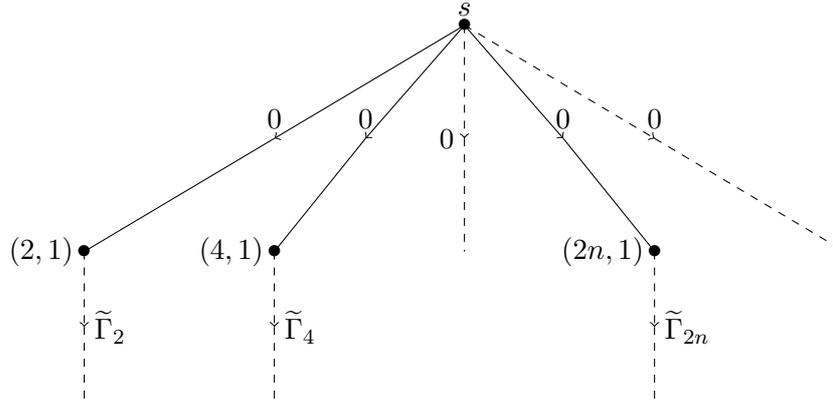
\begin{figure}

\begin{center}

\begin{tikzpicture}[scale=.1]

\draw [->] (50,30) node {$\bullet$} node [anchor=south] {$s$} --
(25,15) node [anchor=south] {$0$};

\draw [-] (25,15) -- (0,0) node {$\bullet$} node [anchor=east]
{$(2,1)$} ;

\draw [->] (50,30) -- (37,15) node [anchor=south] {$0$};

\draw [-] (37,15) -- (25,0) node {$\bullet$} node [anchor=east]
{$(4,1)$} ;

\draw [->][dashed] (50,30) -- (50,15) node [anchor=east] {$0$} ;

\draw [-][dashed] (50,15) -- (50,0) ;

\draw [->] (50,30) -- (63,15) node [anchor=south] {$0$};

\draw [-] (63,15) -- (75,0) node {$\bullet$} node [anchor=east]
{$(2n,1)$}  ;

\draw [->][dashed] (50,30) -- (75,15) node [anchor=south] {$0$} ;

\draw [-][dashed] (75,15) -- (100,0)  ;

\draw [->][dashed] (0,0) -- (0,-10) node [anchor=west]
{$\widetilde{\Gamma}_2$} ;

\draw [-][dashed] (0,-10) -- (0,-20)  ;

\draw [->][dashed] (25,0) -- (25,-10) node [anchor=west]
{$\widetilde{\Gamma}_4$} ;

\draw [-][dashed] (25,-10) -- (25,-20)  ;

\draw [->][dashed] (75,0) -- (75,-10) node [anchor=west]
{$\widetilde{\Gamma}_{2n}$} ;

\draw [-][dashed] (75,-10) -- (75,-20)  ;

%\draw [-] (50,30) -- (30,0) node {$\bullet$} node [anchor=north]
%{$\omega_{2,1}$} ;

\end{tikzpicture}

\caption{The game $\Gamma$ starting from  state $s$}

\end{center}

\end{figure}

\begin{figure}

\begin{center}

\begin{tabular}{rcccrcccrcccrcccrcc}
&C&A&&&C&A&&&C&A&&&C&A&&&C&A\\
\cline{2-3} \cline{6-7} \cline{10-11} \cline{14-15} \cline{18-19}
C&\multicolumn{1}{|c|}{$\underrightarrow{1}$}&\multicolumn{1}{|c|}{0*}&&
C&\multicolumn{1}{|c|}{$\underrightarrow{1}$}&\multicolumn{1}{|c|}{$x_{2n,2}$*}&&
C&\multicolumn{1}{|c|}{$\underrightarrow{1}$}&\multicolumn{1}{|c|}{$x_{2n,m}$*}&&
C&\multicolumn{1}{|c|}{$\underrightarrow{1}$}&\multicolumn{1}{|c|}{$x_{2n,n}$*}&&
C&\multicolumn{1}{|c|}{-1*}&\multicolumn{1}{|c|}{-1*}\\
\cline{2-3} \cline{6-7} \cline{10-11} \cline{14-15} \cline{18-19}
A&\multicolumn{1}{|c|}{0*}&\multicolumn{1}{|c|}{0*}&&
A&\multicolumn{1}{|c|}{$x_{2n,2}$*}&\multicolumn{1}{|c|}{$x_{2n,2}$*}&&
A&\multicolumn{1}{|c|}{$x_{2n,m}$*}&\multicolumn{1}{|c|}{$x_{2n,m}$*}&&
A&\multicolumn{1}{|c|}{$x_{2n,n}$*}&\multicolumn{1}{|c|}{$x_{2n,n}$*}&&
A&\multicolumn{1}{|c|}{-1*}&\multicolumn{1}{|c|}{-1*} \\
\cline{2-3} \cline{6-7} \cline{10-11} \cline{14-15} \cline{18-19}\\
&\multicolumn{2}{c}{$(2n,1)$}&
\multicolumn{1}{c}{$\vphantom{\cdot}$}&
&\multicolumn{2}{c}{$(2n,2)$}& \multicolumn{1}{c}{$\cdots$}&
&\multicolumn{2}{c}{$(2n,m)$}& \multicolumn{1}{c}{$\cdots$}&
&\multicolumn{2}{c}{$(2n,n)$}&
\multicolumn{1}{c}{$\vphantom{\cdot}$}& &\multicolumn{2}{c}{$-1^*$}
\end{tabular}
\end{center}
\caption{The subgame $\widetilde{\Gamma}_{2n}$ starting from
state $(2n,1)$}
\end{figure}

 Notice that in the $2n + 1$ stage game, after a move of player 1 to $\tilde \Gamma_{2n}$, any play is compatible with optimal strategies, in particular those  leading to the sequence of payoffs $2n$ times 0 or $n$ times 1 then $n$ times $-1$.

\subsection{Conjectures.} { }  \\
A natural conjecture is that in any regular game (i.e. where $v_n$ converges uniformly to $v$): \\
for any $ \varepsilon > 0$, there exists $n_0$,  such that for all $n\geq n_0$,  for any initial state $s$, there exists a couple  $(\sigma_n, \tau_n)$ of  $\varepsilon$-optimal strategies in  $G_n (s)$  such that   for any $t\in [0,1]$:
\begin{equation}\label{jj}
3\varepsilon\geq  \frac {1}{n}{\bf E}_{\sigma_n, \tau_n} ^{s}  (  \sum_{m=1}^{[tn]} f_m ) - tv(s) \geq -3 \varepsilon.
\end{equation}
where $[tn]$ stands for the integer part of $tn$ and $f_m$ is  the payoff at stage $m$.\\
A more elaborate conjecture would rely on the existence of an asymptotic game $\Gamma^*$ played in continuous time on $[0,1]$ with value $v$ (as in Section 5.1). We use the representation of the repeated game as a stochastic game trough the recursive structure as above,  see Mertens, Sorin, Zamir \cite{MSZ}, Chapter IV.
The condition  is now the existence of a couple of  strategies  $(\sigma, \tau)$  in the asymptotic game that would depend only  on  the time $t\in [0,1]$ and on  the current state $s$
such that for any $\varepsilon > 0$, there exists $\eta$ with the following property: in any repeated game  where the  (relative) weight of stage $m$ is  $\alpha_m$, with $\{\alpha_m\}$ decreasing and   less than $\eta$,   thus defining a  partition $\Pi$  of $[0,1]$,  the
strategies  $(\sigma_{\Pi}, \tau_{\Pi})$  induced  in the repeated game  by $(\sigma, \tau)$  satisfies $(\ref{jj})$.\\

{\bf Acknowledgment}:\\
This work was done while the three authors were members of the Equipe Combinatoire et Optimisation.\\
Sorin's research was supported by grant ANR-08-BLAN-0294-01 (France).

\end{document}